\documentclass[a4paper,11pt]{amsart}

\usepackage{graphicx}
\usepackage{mathptmx}
\usepackage{amsmath}
\usepackage{amssymb}
\usepackage{enumitem}
\usepackage{xcolor}
\usepackage{pgfplots}

\newmuskip\pFqmuskip

\newcommand*\pFq[6][8]{%
  \begingroup 
  \pFqmuskip=#1mu\relax
  \mathcode`=\string"8000
  \begingroup\lccode`\~=`\,
  \lowercase{\endgroup\let~}\pFqcomma
  F^{#2}_{#3}{\left(\genfrac..{0pt}{}{#4}{#5}\bigg|#6\right)}%
  \endgroup
}
\newcommand{\pFqcomma}{\mskip\pFqmuskip}

\newtheorem{theorem}{Theorem}[section]

\newtheorem{remark}[theorem]{Remark}

\begin{document}

\title[Recurrence relations for fully degenerate Bell polynomials]{A Generalized Recurrence for fully degenerate Bell polynomials}

\author{Taekyun  Kim}
\address{Department of Mathematics, Kwangwoon University, Seoul 139-701, Republic of Korea}
\email{tkkim@kw.ac.kr}
\author{Dae San  Kim}
\address{Department of Mathematics, Sogang University, Seoul 121-742, Republic of Korea}
\email{dskim@sogang.ac.kr}


\begin{abstract}
This paper addresses the unnatural appearance of the two-variable degenerate Fubini polynomials in a recently derived Spivey-type recurrence relation for the fully degenerate Bell polynomials $\mathrm{Bel}_{n,\lambda}(x)$. To solve this, we introduce a new family of polynomials, denoted $B_{n,\lambda}(x)$, which we also call the fully degenerate Bell polynomials, along with their two-variable counterparts, $B_{n,\lambda}(x,y)$. Our main contribution is the derivation of natural Spivey-type recurrence relations for $B_{n,\lambda}(x)$ and $B_{n,\lambda}(x,y)$ using operator methods. We extend these results to the $r$-counterparts, the fully degenerate $r$-Bell polynomials $B_{n,\lambda}^{(r)}(x)$ and $B_{n,\lambda}^{(r)}(x,y)$, providing Dobinski-like, finite sum, operator expressions, and Spivey-type recurrence relations for all the new polynomials.
\end{abstract}

\maketitle

\section{Introduction}
For any nonzero $\lambda\in\mathbb{R}$, the degenerate exponentials are defined by
\begin{equation*}
e_{\lambda}^{x}(t)=\sum_{n=0}^{\infty}(x)_{n,\lambda}\frac{t^{n}}{n!},\quad e_{\lambda}(t)=e_{\lambda}^{1}(t),\quad (\mathrm{see}\ [14-17]),\label{1}
\end{equation*}
where the degenerate falling sequence is given by
$$ (x)_{0,\lambda} = 1,  (x)_{n,\lambda} = x(x - \lambda) \cdots (x - (n-1)\lambda),\quad  (n \geq 1),$$
and satisfies
\begin{equation}\label{1}
(x+y)_{n,\lambda}=\sum_{l=0}^{n}\binom{n}{l}(x)_{l,\lambda}(y)_{n-l,\lambda},\quad(n \ge 0).
\end{equation}
Note that $\displaystyle\lim_{\lambda\rightarrow 0}e_{\lambda}^{x}(t)=e^{xt}\displaystyle$. \par
The degenerate Stirling numbers of the second kind ${n \brace k}_{\lambda}$ are defined by (see \ [10,11,14])
\begin{equation} \label{2}
(x)_{n,\lambda}=\sum_{k=0}^{n}{n \brace k}_{\lambda}(x)_{k},\quad (n\ge 0),\quad  \frac{1}{k!}\big(e_{\lambda}(t)-1 \big)^{k}=\sum_{n=k}^{\infty}{n \brace k}_{\lambda}\frac{t^{n}}{n!}, \quad (k \ge 0).
\end{equation}
For $r\in\mathbb{Z}$ with $r\ge 0$, the degenerate $r$-Stirling numbers of the second kind are given as (see \ [15,17])
\begin{align}
&(x+r)_{n,\lambda}=\sum_{k=0}^{n}{n+r \brace k+r}_{r,\lambda},\ (n\ge 0), \label{3}\\
& \frac{1}{k!}\big(e_{\lambda}(t)-1\big)^{k}e_{\lambda}^{r}(t)=\sum_{n=k}^{\infty}{n+r \brace k+r}_{r,\lambda}\frac{t^{n}}{n!}, \quad (k \ge 0). \nonumber
\end{align}
Note that $\lim_{\lambda\rightarrow 0}{n \brace k}_{\lambda}={n \brace k}$, where ${n \brace k}$ are the Stirling numbers of the second kind defined by
\begin{equation*}
x^{n}=\sum_{k=0}^{n}{n \brace k}(x)_{k},\quad (n\ge 0),\quad (\mathrm{see}\ [4,20]).
\end{equation*}
The Bell polynomials are given by
\begin{equation}
e^{x(e^{t}-1)}=\sum_{n=0}^{\infty}\phi_{n}(x)\frac{t^{n}}{n!},\quad (\mathrm{see}\ [4,5,20]).\label{4}
\end{equation}
By \eqref{4}, we get
\begin{equation}
\phi_{n}(x)=\sum_{k=0}^{n}{n \brace k}x^{k},\quad (\mathrm{see}\ [4,5,20]). \label{5}	
\end{equation}
When $x=1,\ \phi_{n}=\phi_{n}(1)$ are called the Bell numbers. \par
In 2008, Spivey discovered a recurrence relation for Bell numbers:
\begin{equation}
\phi_{n+m}=\sum_{k=0}^{n}\sum_{j=0}^{m}j^{n-k}{m \brace j}\binom{n}{k}\phi_{k},\ (n,m\ge 0),\quad (\mathrm{see}\ [23]). \label{6}	
\end{equation}
It is known that the degenerate Bell polynomials are defined by (see [8,10,21])
\begin{equation}
e^{x(e_{\lambda}(t)-1)}=\sum_{n=0}^{\infty}\phi_{n,\lambda}(x)\frac{t^{n}}{n!},\quad (\mathrm{see}\ [10,11]).\label{7}
\end{equation}
When $x=1,\ \phi_{n,\lambda}=\phi_{n,\lambda}(1)$ are called the degenerate Bell numbers. \par
From \eqref{2} and \eqref{7}, we note that
\begin{equation*}
\phi_{n,\lambda}(x)=\sum_{k=0}^{n}{n \brace k}_{\lambda}x^{k},\ (n\ge 0),\ (\mathrm{see}\ [10,11]).
\end{equation*}
In [10], we obtained a recurrence relation for degenerate Bell polynomials given by
\begin{equation}
\phi_{n+m,\lambda}(x)=\sum_{k=0}^{m}\sum_{l=0}^{n}\binom{n}{l}{m \brace k}_{\lambda}(k-m\lambda)_{n-l,\lambda}x^{k}\phi_{l,\lambda}(x). \label{8}	
\end{equation}
In [12], we considered the fully degenerate Bell polynomials defined by
\begin{equation}
e_{\lambda}\Big(x\big(e_{\lambda}(t)-1\big)\Big)=\sum_{n=0}^{\infty}\mathrm{Bel}_{n,\lambda}(x)\frac{t^{n}}{n!}.\label{9}
\end{equation}
When $x=1,\ \mathrm{Bel}_{n,\lambda}=\mathrm{Bel}_{n,\lambda}(1)$ are called the fully degenerate Bell numbers. We observe here that
\begin{displaymath}
\lim_{\lambda\rightarrow 0}\mathrm{Bel}_{n,\lambda}(x)=\phi_{n}(x),\ (n\ge 0).
\end{displaymath}
Recently, we showed the following Spivey's formula for fully degenerate Bell polynomials (see [12]):
\begin{equation}
\mathrm{Bel}_{n+m,\lambda}(x)=\sum_{k=0}^{m}\sum_{l=0}^{n}(1)_{k,\lambda}{m \brace k}_{\lambda}\binom{n}{l}x^{k}F_{n-l,\lambda}^{(k)}(-\lambda l, k-m\lambda)\mathrm{Bel}_{l,\lambda}(x), \label{10}	
\end{equation}
where $F_{n,\lambda}^{(k)}(x,y)$ are the two variable degenerate Fubini polynomials of order $k$ given by
\begin{displaymath}
\bigg(\frac{1}{1-x(e_{\lambda}(t)-1)}\bigg)^{k}e_{\lambda}^{y}(t)=\sum_{n=0}^{\infty}F_{n,\lambda}^{(k)}(x,y)\frac{t^{n}}{n!},\quad (\mathrm{see}\ [12]).
\end{displaymath}
Note that
\begin{displaymath}
\lim_{\lambda\rightarrow 0}F_{n,\lambda}^{(1)}(x,0)=F_{n}(x),
\end{displaymath}
where $F_{n}(x)$ are the Fubini polynomials defined by
\begin{displaymath}
F_{n}(x)=\sum_{k=0}^{n}{n \brace k}k!x^{k},\quad (n\ge 0),\quad (\mathrm{see}\ [12]).
\end{displaymath}
As general references of this paper, the reader may refer to [4,13,18-20,22]. \par
\vspace{0.1in}
The study of degenerate versions of special numbers and polynomials traces back to the pioneering work of Carlitz [3]. This field has recently experienced a resurgence of interest, with numerous degenerate counterparts of various special numbers and polynomials being investigated through diverse methods (see [1, 9, 11, 12, 14–17] and their references).
A Spivey-type recurrence relation was recently established for the fully degenerate Bell polynomials $Bel_{n,\lambda}(x)$ (see \eqref{10}). However, this particular recurrence relation \eqref{10} is somewhat unnatural (see \eqref{6}, \eqref{8}) because it involves the two-variable degenerate Fubini polynomials of order $k$, which are unrelated to the fully degenerate Bell polynomials. To resolve this issue, we introduce a new set of polynomials, which we also term fully degenerate Bell polynomials but denote with the distinct notation $B_{n, \lambda}(x)$. Our primary goal is to derive Spivey-type recurrence relations for these new polynomials, $B_{n,\lambda}(x)$, and their two-variable counterparts, $B_{n,\lambda}(x,y)$. Furthermore, we extend our investigation to find Spivey-type recurrence relations for their $r$-counterparts: the fully degenerate $r$-Bell polynomials $B_{n,\lambda}^{(r)}(x)$, and the two-variable fully degenerate r-Bell polynomials, $B_{n,\lambda}^{(r)}(x,y)$.
Our derivation utilizes the `multiplication by $x$' operator $X$ and the differentiation operator $D=\frac{d}{dx}$, which satisfy the fundamental commutation relation $DX-XD=1$. Alternatively, one could employ Boson operators (see [2, 15, 16]). For prior results on Spivey-type relations, the reader may refer to [5, 11, 12, 23]. Specifically, for $q$-analogues of Spivey’s Bell number formula, we direct the reader to [6, 7]. \par
This paper is structured to clearly present our definitions and main findings regarding the newly introduced fully degenerate Bell polynomials. In Section 1, we begin by recalling essential facts and notation that will be used throughout the paper. Section 2 is the core of the paper, where we introduce the new fully degenerate Bell polynomials $B_{n,\lambda}(x)$ and the two-variable fully degenerate Bell polynomials $B_{n,\lambda}(x,y)$.
Our results for these polynomials include: \par
\vspace{0.1in}
$\bullet$  $B_{n,\lambda}(x)$ (Theorems 2.1, 2.3, 2.4): We provide a finite sum expression, a Dobinski-like formula, and an operator expression.\par
$\bullet$  $B_{n,\lambda}(x,y)$ (Theorems 2.5, 2.6, 2.7): Similarly, we establish a finite sum expression, a Dobinski-like formula, and an operator expression.\par
\vspace{0.1in}
In Theorem 2.8, we derive the central Spivey-type recurrence relations for $B_{n,\lambda}(x,y)$ and $B_{n,\lambda}(x)$: \\
\begin{align*}
&B_{n+m,\lambda}(x,y)= \sum_{j=0}^{m}\sum_{k=0}^{n}{m\brace j}_{\lambda}\binom{n}{k}(j-m\lambda)_{n-k,\lambda}(y)_{j,\lambda}x^{j}B_{k,\lambda}(x,y-j\lambda), \\
&B_{n+m,\lambda}(x)= \sum_{j=0}^{m}\sum_{k=0}^{n}{m\brace j}_{\lambda}\binom{n}{k}(j-m\lambda)_{n-k,\lambda}(1)_{j,\lambda}x^{j}B_{k,\lambda}(x,1-j\lambda).
\end{align*}
\vspace{0.1in}
Next, we define the fully degenerate $r$-Bell polynomials $B_{n,\lambda}^{(r)}(x)$
and the two-variable fully degenerate $r$-Bell polynomials $B_{n,\lambda}^{(r)}(x,y)$. We provide analogous results for these polynomials: \par
\vspace{0.1in}
$\bullet$ $B_{n,\lambda}^{(r)}(x)$ (Theorems 2.9, 2.10, 2.11): A finite sum expression, a Dobinski-like formula, and an operator expression. \par
$\bullet$ $B_{n,\lambda}^{(r)}(x,y)$ (Theorems 2.12, 2.13, 2.14): A finite sum expression, a Dobinski-like formula, and an operator expression. \par
\vspace{0.1in}
Finally, in Theorem 2.15, we present the Spivey-type recurrence relations for $B_{n,\lambda}^{(r)}(x,y)$ and $B_{n,\lambda}^{(r)}(x)$: \\
\begin{align*}
&B_{n+m,\lambda}^{(r)}(x,y)=\sum_{k=0}^{m}\sum_{l=0}^{n}\binom{n}{l}{m+r \brace k+r}_{r,\lambda}(k-m\lambda)_{n-l,\lambda}(y)_{k,\lambda}x^{k}B_{l,\lambda}^{(r)}(x,y-k\lambda), \\
&B_{n+m,\lambda}^{(r)}(x)=\sum_{k=0}^{m}\sum_{l=0}^{n}\binom{n}{l}{m+r \brace k+r}_{r,\lambda}(k-m\lambda)_{n-l,\lambda}(1)_{k,\lambda}x^{k}B_{l,\lambda}^{(r)}(x,1-k\lambda).
\end{align*} \par
In Section 3, we conclude the paper by summarizing the key findings and suggesting avenues for future research.

\section{Recurrence relations for fully degenerate Bell polynomials}
In this section, we use the operators $X$ and $D$ defined by $Xf(x)=xf(x)$ and $Df(x)=\frac{d}{dx}f(x)$. They satisfy the commutation relation $DX-XD=1$. \par
Now, we redefine the polynomials $\mathrm{Bel}_{n,\lambda}(x)$ in \eqref{9}, again called the {\it{fully degenerate Bell polynomials}} but denoted by the different notation $B_{n,\lambda}(x)$, by
\begin{equation}
\frac{1}{e_{\lambda}(x)}e_{\lambda}\big(xe_{\lambda}(t)\big)=\sum_{n=0}^{\infty}B_{n,\lambda}(x)\frac{t^{n}}{n!}. \label{11}	
\end{equation}
Note from \eqref{4} that
\begin{align}
\lim_{\lambda\rightarrow 0}\sum_{n=0}^{\infty}B_{n,\lambda}(x)\frac{t^{n}}{n!}&=\lim_{\lambda\rightarrow 0}e_{\lambda}^{-1}(x)e_{\lambda}\big(xe_{\lambda}(t)\big)=e^{x(e^{t}-1)}\label{12}\\
&=\sum_{n=0}^{\infty}\phi_{n}(x)\frac{t^{n}}{n!}.\nonumber
\end{align}
Thus, by \eqref{12}, we get
\begin{displaymath}
\lim_{\lambda\rightarrow 0}B_{n,\lambda}(x)=\phi_{n}(x),\ (n\ge 0).
\end{displaymath}
When $x=1,\ B_{n,\lambda}=B_{n,\lambda}(1)$ are called the fully degenerate Bell numbers. \par
From \eqref{2} and \eqref{11}, we note that
\begin{align}
&\sum_{n=0}^{\infty}B_{n,\lambda}(x)\frac{t^{n}}{n!}=\frac{1}{e_{\lambda}(x)}e_{\lambda}\big(xe_{\lambda}(t)\big)= \bigg(\frac{1+\lambda x+\lambda x(e_{\lambda}(t)-1)}{1+\lambda x}\bigg)^{\frac{1}{\lambda}}\label{13}\\
&=\bigg(1+\lambda\frac{x}{1+\lambda x}\big(e_{\lambda}(t)-1\big)\bigg)^{\frac{1}{\lambda}}=e_{\lambda}\bigg(\frac{x}{1+\lambda x}(e_{\lambda}(t)-1)\bigg)\nonumber\\
&=\sum_{k=0}^{\infty}(1)_{k,\lambda}\bigg(\frac{x}{1+\lambda x}\bigg)^{k}\frac{1}{k!}\big(e_{\lambda}(t)-1\big)^{k}=\sum_{k=0}^{\infty}(1)_{k,\lambda}\bigg(\frac{x}{1+\lambda x}\bigg)^{k}\sum_{n=k}^{\infty}{n \brace k}_{\lambda}\frac{t^{n}}{n!} \nonumber\\
&=\sum_{n=0}^{\infty}\sum_{k=0}^{n}(1)_{k,\lambda}\bigg(\frac{x}{1+\lambda x}\bigg)^{k}{n \brace k}_{\lambda}\frac{t^{n}}{n!}. \nonumber
\end{align}
Therefore, by comparing the coefficients on both sides of \eqref{13}, we obtain the following theorem.
\begin{theorem}
For $n\ge 0$, we have
\begin{equation}
B_{n,\lambda}(x)=\sum_{k=0}^{n}(1)_{k,\lambda}\bigg(\frac{x}{1+\lambda x}\bigg)^{k}{n \brace k}_{\lambda}. \label{14}
\end{equation}
In particular, for $x=1$, we get
\begin{displaymath}
B_{n,\lambda}=\sum_{k=0}^{n}(1)_{k,\lambda}\bigg(\frac{1}{1+\lambda}\bigg)^{k}{n \brace k}_{\lambda}.
\end{displaymath}
\end{theorem}
\begin{remark}
From \eqref{9}, we note that $\mathrm{Bel}_{n,\lambda}(x)$ is given by
\begin{equation*}
\mathrm{Bel}_{n,\lambda}(x)=\sum_{k=0}^{n}(1)_{k,\lambda}x^{k}{n \brace k}_{\lambda}.
\end{equation*}
\end{remark}
By \eqref{11} and Taylor expansion, we get
\begin{align}
\sum_{n=0}^{\infty}B_{n,\lambda}(x)\frac{t^{n}}{n!}&=\frac{1}{e_{\lambda}(x)}e_{\lambda}\big(xe_{\lambda}(t)\big)=\frac{1}{e_{\lambda}(x)}\sum_{k=0}^{\infty}\frac{(1)_{k,\lambda}x^{k}}{k!}e_{\lambda}^{k}(t) \label{15}\\
&=\frac{1}{e_{\lambda}(x)}\sum_{k=0}^{\infty}\frac{(1)_{k,\lambda}}{k!}x^{k}\sum_{n=0}^{\infty}(k)_{n,\lambda}\frac{t^{n}}{n!}\nonumber\\
&=\sum_{n=0}^{\infty}\frac{1}{e_{\lambda}(x)}\sum_{k=0}^{\infty}\frac{(1)_{k,\lambda}}{k!}(k)_{n,\lambda}x^{k}\frac{t^{n}}{n!}. \nonumber	
\end{align}
Thus, by \eqref{15}, we obtain a Dobinski-like formula for $B_{n,\lambda}(x)$.
\begin{theorem}
For $n\ge 0$, we have
\begin{equation}
B_{n,\lambda}(x)=\frac{1}{e_{\lambda}(x)}\sum_{k=0}^{\infty}(1)_{k,\lambda}\frac{(k)_{n,\lambda}}{k!}x^{k}.\label{16}
\end{equation}
In particular, for $x=1$, we get
\begin{displaymath}
B_{n,\lambda}=\frac{1}{e_{\lambda}(1)}\sum_{k=0}^{\infty}(1)_{k,\lambda}\frac{(k)_{n,\lambda}}{k!}.
\end{displaymath}
\end{theorem}
Note that letting $\lambda \rightarrow 0$ in \eqref{16} gives Dobinski's formula:
\begin{displaymath}
\phi_{n}(x)=\frac{1}{e^{x}}\sum_{k=0}^{\infty}\frac{k^{n}}{k!}x^{k},\ (n\ge 0).
\end{displaymath}
Using \eqref{2} and $DX-XD=1$, we show that
\begin{equation}
(XD)_{n,\lambda}=\sum_{k=0}^{n}{n \brace k}_{\lambda}X^{k}D^{k},\ (n\ge 0),\ (\mathrm{see}\ [15,16]).\label{17}
\end{equation}
From \eqref{16}, we note that
\begin{align}
\frac{1}{e_{\lambda}(x)}(XD)_{n,\lambda}e_{\lambda}(x)&=\frac{1}{e_{\lambda}(x)}\sum_{k=0}^{\infty}\frac{(1)_{k,\lambda}}{k!}(XD)_{n,\lambda}x^{k}\label{18}\\
&=\frac{1}{e_{\lambda}(x)}\sum_{k=0}^{\infty}(1)_{k,\lambda}\frac{(k)_{n,\lambda}}{k!}x^{k}\nonumber\\
&=B_{n,\lambda}(x),\quad (n\ge 0). \nonumber
\end{align}
Therefore, by \eqref{18}, we obtain the following theorem.
\begin{theorem}
For $n\ge 0$, in terms of operators we have
\begin{equation*}
B_{n,\lambda}(x)=\frac{1}{e_{\lambda}(x)}(XD)_{n,\lambda}e_{\lambda}(x). 	
\end{equation*}
\end{theorem}
Now, we observe that
\begin{equation}
(XD)_{n+m,\lambda}=\big(XD-m\lambda\big)_{n,\lambda}(XD)_{m,\lambda},\quad (m,n\ge 0). \label{19}
\end{equation}
For $k\in\mathbb{N}$, we observe that
\begin{align}
DX^{k}-X^{k}D&=\big(DX^{k-1}-X^{k-1}D)X+X^{k-1}(DX-XD)  \label{20}\\
&=\big(DX^{k-1}-X^{k-1}D\big)X+X^{k-1}\nonumber\\
&=\big(DX^{k-2}-X^{k-2}D\big)X^{2}+2X^{k-1}=\cdots \nonumber\\
&=\big(DX-XD\big)X^{k-1}+(k-1)X^{k-1}=kX^{k-1}. \nonumber
\end{align}
By \eqref{20}, we get
\begin{align}
\big(XD\big)X^{k}=X\big(DX^{k}\big)=X^{k}\big(XD+k\big),\quad  (k\in\mathbb{N}). \label{21}
\end{align}
For $m,n\ge 0$, by \eqref{21}, we get
\begin{align}
\big(XD-m\lambda\big)_{n,\lambda}X^{j}&=\big(XD-m\lambda\big)\big(XD-m\lambda-\lambda\big)\cdots\big(XD-m\lambda-(n-1)\lambda\big)X^{j}\label{22} \\
&=X^{j}\big(XD+j-m\lambda\big)_{n,\lambda}=X^{j}\sum_{k=0}^{n}\binom{n}{k}(XD)_{k,\lambda}(j-m\lambda)_{n-k,\lambda}. \nonumber
\end{align}
From \eqref{1}, \eqref{17}, \eqref{19}, \eqref{21} and \eqref{22}, we have
\begin{align}
\big(XD\big)_{n+m,\lambda}&=\big(XD-m\lambda\big)_{n,\lambda}\big(XD\big)_{m,\lambda}=\sum_{j=0}^{m}{m \brace j}_{\lambda}\big(XD-m\lambda\big)_{n,\lambda}X^{j}D^{j} \label{23}\\
&=\sum_{j=0}^{m}{m \brace j}_{\lambda}X^{j}\big(XD+j-m\lambda\big)_{n,\lambda}D^{j}\nonumber\\
&=\sum_{j=0}^{m}\sum_{k=0}^{n}{m \brace j}_{\lambda}\binom{n}{k}(j-m\lambda)_{n-k,\lambda}X^{j}(XD)_{k,\lambda}D^{j}.\nonumber
\end{align}

Now, we define the {\it{two variable fully degenerate Bell polynomials}} by
\begin{equation}
\frac{1}{e_{\lambda}(x)}e_{\lambda}^{y}\big(xe_{\lambda}(t)\big)=\sum_{n=0}^{\infty}B_{n,\lambda}(x,y)\frac{t^{n}}{n!}.\label{24}	
\end{equation}
From \eqref{24} and working just as in \eqref{13}, we have
\begin{align}
\sum_{n=0}^{\infty}B_{n,\lambda}(x,y)\frac{t^{n}}{n!}&=\frac{1}{e_{\lambda}(x)}e_{\lambda}^{y}\big(xe_{\lambda}(t)\big)\label{25}\\
&=e_{\lambda}^{y-1}(x)e_{\lambda}^{y}\bigg(\frac{x}{1+\lambda x}\big(e_{\lambda}(t)-1\big)\bigg)\nonumber\\
&=e_{\lambda}^{y-1}(x)\sum_{k=0}^{\infty}(y)_{k,\lambda}\bigg(\frac{x}{1+\lambda x}\bigg)^{k}\frac{1}{k!}\big(e_{\lambda}(t)-1\big)^{k} \nonumber \\
&=e_{\lambda}^{y-1}(x)\sum_{k=0}^{\infty}(y)_{k,\lambda}\bigg(\frac{x}{1+\lambda x}\bigg)^{k}\sum_{n=k}^{\infty}{n \brace k}_{\lambda}\frac{t^{n}}{n!} \nonumber \\
&=e_{\lambda}^{y-1}(x)\sum_{n=0}^{\infty}\sum_{k=0}^{n}{n \brace k}_{\lambda}(y)_{k,\lambda}\bigg(\frac{x}{1+\lambda x}\bigg)^{k}\frac{t^{n}}{n!}. \nonumber
\end{align}
Therefore, by \eqref{14} and \eqref{25}, we obtain the following theorem.
\begin{theorem}
For $n\ge 0$, we have
\begin{equation*}
B_{n,\lambda}(x,y)=e_{\lambda}^{y-1}(x)\sum_{k=0}^{n}{n \brace k}_{\lambda}(y)_{k,\lambda}\bigg(\frac{x}{1+\lambda x}\bigg)^{k}. 	
\end{equation*}
When $y=1$, we get
\begin{equation}
B_{n,\lambda}(x,1)=\sum_{k=0}^{n}{n \brace k}_{\lambda}(1)_{k,\lambda}\bigg(\frac{x}{1+\lambda x}\bigg)^{k}=B_{n,\lambda}(x). \label{26}
\end{equation}
\end{theorem}
By \eqref{24} and Taylor expansion, we get
\begin{align}
\sum_{n=0}^{\infty}B_{n,\lambda}(x,y)\frac{t^{n}}{n!}&=\frac{1}{e_{\lambda}(x)}e_{\lambda}^{y}\big(xe_{\lambda}(t)\big) \label{27} \\
&=\frac{1}{e_{\lambda}(x)}\sum_{k=0}^{\infty}\frac{(y)_{k,\lambda}}{k!}x^{k}e_{\lambda}^{k}(t)\nonumber\\
&=\sum_{n=0}^{\infty}\frac{1}{e_{\lambda}(x)}\sum_{k=0}^{\infty}\frac{(y)_{k,\lambda}}{k!}(k)_{n,\lambda}x^{k}\frac{t^{n}}{n!}. \nonumber	
\end{align}
Therefore, by \eqref{27}, we obtain a Dobinski-like formula for $B_{n,\lambda}(x,y)$.
\begin{theorem}
For $n\ge 0$, we have
\begin{equation}
B_{n,\lambda}(x,y)=\frac{1}{e_{\lambda}(x)}\sum_{k=0}^{\infty}\frac{(y)_{k,\lambda}}{k!}(k)_{n,\lambda}x^{k}.\label{28}
\end{equation}
\end{theorem}
From \eqref{28}, we have
\begin{align}
\frac{1}{e_{\lambda}(x)}\big(XD\big)_{n,\lambda}e_{\lambda}^{y}(x)&=\frac{1}{e_{\lambda}(x)}(XD)_{n,\lambda}\sum_{k=0}^{\infty}\frac{(y)_{k,\lambda}}{k!}x^{k} \label{29} \\
&=\frac{1}{e_{\lambda}(x)}\sum_{k=0}^{\infty}\frac{(y)_{k,\lambda}}{k!}(k)_{n,\lambda}x^{k}\nonumber\\
&=B_{n,\lambda}(x,y). \nonumber
\end{align}
Thus, by \eqref{29}, we have the following theorem.
\begin{theorem}
For $n \ge 0$, in terms of operators we have
\begin{equation}
B_{n,\lambda}(x,y)=\frac{1}{e_{\lambda}(x)}\big(XD\big)_{n,\lambda}e_{\lambda}^{y}(x)\label{30}
\end{equation}
\end{theorem}
By \eqref{23}, \eqref{29} and \eqref{30}, and noting that $D^{j}e_{\lambda}^{y}(x)=(y)_{j,\lambda}e_{\lambda}^{y-j \lambda}(x)$, we get
\begin{align}
B_{n+m,\lambda}(x,y)&=\frac{1}{e_{\lambda}(x)}(XD)_{n+m,\lambda}e_{\lambda}^{y}(x) \label{31}\\
&=\frac{1}{e_{\lambda}(x)}\sum_{j=0}^{m}\sum_{k=0}^{n}{m \brace j}_{\lambda}\binom{n}{k}(j-m\lambda)_{n-k,\lambda}X^{j}(XD)_{k,\lambda}D^{j}e_{\lambda}^{y}(x) \nonumber\\
&=\sum_{j=0}^{m}\sum_{k=0}^{n}{m \brace j}_{\lambda}\binom{n}{k}(j-m\lambda)_{n-k,\lambda}(y)_{j,\lambda}X^{j}\frac{1}{{e_{\lambda}(x)}}(XD)_{k,\lambda}e_{\lambda}^{y-j\lambda }(x)\nonumber \\
&=\sum_{j=0}^{m}\sum_{k=0}^{n}{m \brace j}_{\lambda}\binom{n}{k}(j-m\lambda)_{n-k,\lambda}(y)_{j,\lambda}x^{j}B_{k,\lambda}(x,y-j \lambda).\nonumber
\end{align}

Therefore, by \eqref{26} and \eqref{31}, we obtain Spivey-type recurrence relations for two variable fully degenerate Bell polynomials and fully degenerate Bell polynomials.
\begin{theorem}
For $n,m\ge 0$, we have
\begin{equation*}
B_{n+m,\lambda}(x,y)= \sum_{j=0}^{m}\sum_{k=0}^{n}{m\brace j}_{\lambda}\binom{n}{k}(j-m\lambda)_{n-k,\lambda}(y)_{j,\lambda}x^{j}B_{k,\lambda}(x,y-j\lambda), \end{equation*}
and
\begin{equation*}
B_{n+m,\lambda}(x)= \sum_{j=0}^{m}\sum_{k=0}^{n}{m\brace j}_{\lambda}\binom{n}{k}(j-m\lambda)_{n-k,\lambda}(1)_{j,\lambda}x^{j}B_{k,\lambda}(x,1-j\lambda).
\end{equation*}
\end{theorem}
Letting $\lambda \rightarrow 0$, we have (see \eqref{5}, \eqref{14}, \eqref{24}, \eqref{32})
\begin{align*}
&\phi_{n+m}(x,y)=\sum_{j=0}^{m}\sum_{k=0}^{n}{m \brace j}\binom{n}{k}j^{n-k}y^{j}x^{j}\phi_{k}(x,y), \\
&\phi_{n+m}(x)=\sum_{j=0}^{m}\sum_{k=0}^{n}{m \brace j}\binom{n}{k}j^{n-k}x^{j}\phi_{k}(x),
\end{align*}
where the two variable Bell polynomials $\phi_{n}(x,y)$ are defined by
\begin{equation}
e^{x(ye^{t}-1)}=\sum_{n=0}^{\infty}\phi_{n}(x,y)\frac{t^{n}}{n!}.\label{32}
\end{equation}
For $r\ge 0$, we define the {\it{fully degenerate $r$-Bell polynomials}} by
\begin{equation}
e_{\lambda}^{-1}(x)e_{\lambda}\big(xe_{\lambda}(t)\big)e_{\lambda}^{r}(t)=\sum_{n=0}^{\infty}B_{n,\lambda}^{(r)}(x)\frac{t^{n}}{n!}. \label{33}
\end{equation}
When $r=0,\ B_{n,\lambda}(x)=B_{n,\lambda}^{(0)}(x)$ are the fully degenerate Bell polynomials. \par
From \eqref{3}, \eqref{13} and \eqref{33}, we have
\begin{align}
e_{\lambda}^{-1}(x)e_{\lambda}\big(xe_{\lambda}(t)\big)e_{\lambda}^{r}(t)
&=e_{\lambda}\bigg(\frac{x}{1+\lambda x}\big(e_{\lambda}(t)-1\big)\bigg)e_{\lambda}^{r}(t) \label{34}\\
&=\sum_{k=0}^{\infty}(1)_{k,\lambda}\bigg(\frac{x}{1+x\lambda x}\bigg)^{k}\frac{1}{k!}\big(e_{\lambda}(t)-1\big)^{k}e_{\lambda}^{r}(t)\nonumber\\
&=\sum_{k=0}^{\infty}(1)_{k,\lambda}\bigg(\frac{x}{1+\lambda x}\bigg)^{k}\sum_{n=k}^{\infty}{n+r \brace k+r}_{r,\lambda}\frac{t^{n}}{n!}\nonumber\\
&=\sum_{n=0}^{\infty}\sum_{k=0}^{n}(1)_{k,\lambda}\bigg(\frac{x}{1+\lambda x}\bigg)^{k}{n+r \brace k+r}_{r,\lambda}\frac{t^{n}}{n!}. \nonumber
\end{align}
Therefore, by \eqref{33} and \eqref{34}, we obtain the following theorem.
\begin{theorem}
For $n\ge 0$, we have
\begin{equation}
B_{n,\lambda}^{(r)}(x)= \sum_{k=0}^{n}(1)_{k,\lambda}\bigg(\frac{x}{1+\lambda x}\bigg)^{k}{n+r \brace k+r}_{r,\lambda}. \label{35}
\end{equation}
\end{theorem}
From \eqref{33}, we have
\begin{align}
\sum_{n=0}^{\infty}B_{n,\lambda}^{(r)}(x)\frac{t^{n}}{n!}&=\frac{1}{e_{\lambda}(x)}e_{\lambda}\big(xe_{\lambda}(t)\big)e_{\lambda}^{r}(t) \label{36} \\
&=\frac{1}{e_{\lambda}(x)}\sum_{k=0}^{\infty}\frac{(1)_{k,\lambda}}{k!}x^{k}e_{\lambda}^{k+r}(t) \nonumber\\
&=\sum_{n=0}^{\infty}\frac{1}{e_{\lambda}(x)}\sum_{k=0}^{\infty}\frac{(1)_{k,\lambda}}{k!}x^{k}(k+r)_{n,\lambda}\frac{t^{n}}{n!}. \nonumber
\end{align}
Therefore, by \eqref{36}, we obtain a Dobinski-like formula for $B_{n,\lambda}^{(r)}(x)$.
\begin{theorem}
For $n\ge 0$, we have
\begin{equation}
B_{n,\lambda}^{(r)}(x)= \frac{1}{e_{\lambda}(x)}\sum_{k=0}^{\infty}\frac{(1)_{k,\lambda}}{k!}(k+r)_{n,\lambda} x^{k}.\label{37}
\end{equation}
\end{theorem}
For $r\ge 0$, by \eqref{37}, we get
\begin{align}
\frac{1}{e_{\lambda}(x)}\big(XD+r\big)_{n,\lambda}e_{\lambda}(x)&=\frac{1}{e_{\lambda}(x)}\sum_{k=0}^{\infty}\frac{(1)_{k,\lambda}}{k!}\big(XD+r\big)_{n,\lambda}x^{k}\label{38}\\
&=\frac{1}{e_{\lambda}(x)}\sum_{k=0}^{\infty}\frac{(1)_{k,\lambda}}{k!}(k+r)_{n,\lambda}x^{k}\nonumber\\
&=B_{n,\lambda}^{(r)}(x). \nonumber
\end{align}
Therefore, by \eqref{38}, we obtain the following theorem.
\begin{theorem}
For $n\ge 0$, we have
\begin{displaymath}
B_{n,\lambda}^{(r)}(x)=\frac{1}{e_{\lambda}(x)}\big(XD+r\big)_{n,\lambda}e_{\lambda}(x).
\end{displaymath}
\end{theorem}
In view of \eqref{24}, we define the {\it{two variable fully degenerate $r$-Bell polynomials}}
by
\begin{equation}
\frac{1}{e_{\lambda}(x)}e_{\lambda}^{y}\big(xe_{\lambda}(t)\big)e_{\lambda}^{r}(t)=\sum_{n=0}^{\infty}B_{n,\lambda}^{(r)}(x,y)\frac{t^{n}}{n!}. \label{39}
\end{equation}
When $y=1$, we have $B_{n,\lambda}^{(r)}(x,1)=B_{n,\lambda}^{(r)}(x)$. \par
From \eqref{3} and \eqref{39}, we note that
\begin{align}
&\sum_{n=0}^{\infty}B_{n,\lambda}^{(r)}(x,y)\frac{t^{n}}{n!}=\frac{1}{e_{\lambda}(x)}e_{\lambda}^{y}\big(xe_{\lambda}(t)\big)e_{\lambda}^{r}(t) \label{40}\\
&=e_{\lambda}^{y-1}(x)e_{\lambda}^{y}\Big(\frac{x}{1+\lambda x}(e_{\lambda}(t)-1)\Big)e_{\lambda}^{r}(t)\nonumber\\
&=e_{\lambda}^{y-1}(x)\sum_{k=0}^{\infty}(y)_{k,\lambda}\Big(\frac{x}{1+\lambda x}\Big)^{k}\frac{1}{k!}\big(e_{\lambda}(t)-1\big)^{k}e_{\lambda}^{r}(t)\nonumber\\
&=e_{\lambda}^{y-1}(x)\sum_{k=0}^{\infty}(y)_{k,\lambda}\Big(\frac{x}{1+\lambda x}\Big)^{k}\sum_{n=k}^{\infty}{n+r \brace k+r}_{r,\lambda}\frac{t^{n}}{n!}\nonumber\\
&=\sum_{n=0}^{\infty}e_{\lambda}^{y-1}(x)\sum_{k=0}^{n}(y)_{k,\lambda}\Big(\frac{x}{1+\lambda x}\Big)^{k}{n+r \brace k+r}_{r,\lambda}\frac{t^{n}}{n!}. \nonumber
\end{align}
Therefore, by \eqref{35} and \eqref{40}, we obtain the following theorem.
\begin{theorem}
For $n\ge 0$, we have
\begin{displaymath}
B_{n,\lambda}^{(r)}(x,y)= e_{\lambda}^{y-1}(x)\sum_{k=0}^{n}(y)_{k,\lambda}\Big(\frac{x}{1+\lambda x}\Big)^{k}{n+r \brace k+r}_{r,\lambda}.
\end{displaymath}
In particular, for $y=1$, we get
\begin{equation}
B_{n,\lambda}^{(r)}(x,1)=\sum_{k=0}^{n}(1)_{k,\lambda}\Big(\frac{x}{1+\lambda x}\Big)^{k}{n+r \brace k+r}_{r,\lambda}=B_{n,\lambda}^{(r)}(x),\ (n\ge 0). \label{41}
\end{equation}
\end{theorem}
By \eqref{39} and Taylor expansion, we get a Dobinski-like formula for $B_{n,\lambda}^{(r)}(x,y)$:
\begin{align}
\sum_{n=0}^{\infty}B_{n,\lambda}^{(r)}(x,y)\frac{t^{n}}{n!}&=\frac{1}{e_{\lambda}(x)}e_{\lambda}^{y}\big(xe_{\lambda}(t)\big)e_{\lambda}^{r}(t) \label{42}\\
&=\sum_{k=0}^{\infty}\frac{1}{e_{\lambda}(x)}\frac{(y)_{k,\lambda}}{k!}x^{k}e_{\lambda}^{r+k}(t)\nonumber\\
&=\sum_{n=0}^{\infty}\frac{1}{e_{\lambda}(x)}\sum_{k=0}^{\infty}\frac{(y)_{k,\lambda}}{k!}x^{k}(r+k)_{n,\lambda}\frac{t^{n}}{n!}. \nonumber	
\end{align}
From \eqref{42}, we have a Dobinski-like formula for $B_{n,\lambda}^{(r)}(x,y)$.
\begin{theorem}
For $n \ge 0$, we have
\begin{equation}
B_{n,\lambda}^{(r)}(x,y)=\frac{1}{e_{\lambda}(x)}\sum_{k=0}^{\infty}\frac{(y)_{k,\lambda}}{k!}(r+k)_{n,\lambda}x^{k},\quad (n\ge 0). \label{43}	
\end{equation}
\end{theorem}
By \eqref{43}, we get
\begin{align}
\frac{1}{e_{\lambda}(x)}\big(XD+r\big)_{n,\lambda}e_{\lambda}^{y}(x)&=\frac{1}{e_{\lambda}(x)}\sum_{k=0}^{\infty}\frac{(y)_{k,\lambda}}{k!}\big(XD+r\big)_{n,\lambda}x^{k}\label{44}\\
&=\frac{1}{e_{\lambda}(x)}\sum_{k=0}^{\infty}\frac{(y)_{k,\lambda}}{k!}(k+r)_{n,\lambda}x^{k}\nonumber\\
&=B_{n,\lambda}^{(r)}(x,y),\quad (n\ge 0).\nonumber
\end{align}
Therefore, by \eqref{44}, we obtain the following theorem.
\begin{theorem}
For $n\ge 0$, in terms of operators we have
\begin{equation}
B_{n,\lambda}^{(r)}(x,y)=\frac{1}{e_{\lambda}(x)}\big(XD+r\big)_{n,\lambda}e_{\lambda}^{y}(x). \label{45}
\end{equation}
\end{theorem}
For $m,n\ge 0$, by using \eqref{3} and $DX-XD=1$, we have
\begin{align}
\big(XD+r\big)_{n+m,\lambda}&=\big(XD+r-m\lambda\big)_{n,\lambda}\big(XD+r\big)_{m,\lambda} \label{46}\\
&=\big(XD+r-m\lambda\big)_{n,\lambda}\sum_{k=0}^{m}{m+r \brace k+r}_{r,\lambda}X^{k}D^{k}.\nonumber
\end{align}
By \eqref{21}, we get
\begin{align}
\big(XD+r-m\lambda\big)_{n,\lambda}X^{k}&=X^{k}\big(XD+r+k-m\lambda\big)_{n,\lambda} \label{47}\\
&=X^{k}\sum_{l=0}^{n}\binom{n}{l}\big(XD+r\big)_{l,\lambda}\big(k-m\lambda\big)_{n-l,\lambda}.\nonumber
\end{align}
From \eqref{1}, \eqref{46} and \eqref{47}, we note that
\begin{align}
&\big(XD+r\big)_{n+m,\lambda}=\sum_{k=0}^{m}{m+r \brace k+r}_{r,\lambda}\big(XD+r-m\lambda\big)_{n,\lambda}X^{k}D^{k}\label{48}\\
&=\sum_{k=0}^{m}{m+r \brace k+r}_{r,\lambda}X^{k}\big(XD+r+k-m\lambda\big)_{n,\lambda}D^{k}\nonumber\\
&=\sum_{k=0}^{m}\sum_{l=0}^{n}{m+r \brace k+r}_{r,\lambda}\binom{n}{l}(k-m\lambda)_{n-l,\lambda}X^{k}\big(XD+r\big)_{l,\lambda}D^{k}. \nonumber	
\end{align}
By \eqref{45}, and noting that $D^{k}e_{\lambda}^{y}(x)=(y)_{k,\lambda}e_{\lambda}^{y-k \lambda}(x)$, we get
\begin{align}
&B_{n+m,\lambda}^{(r)}(x,y)=\frac{1}{e_{\lambda}(x)}\big(XD+r\big)_{n+m,\lambda}e_{\lambda}^{y}(x)\label{49}\\
&=\sum_{k=0}^{m}\sum_{l=0}^{n}{m+r \brace k+r}_{r,\lambda}\binom{n}{l}(k-m\lambda)_{n-l,\lambda}\frac{1}{e_{\lambda}(x)}X^{k}\big(XD+r\big)_{l,\lambda}D^{k}e_{\lambda}^{y}(x)\nonumber\\
&=\sum_{k=0}^{m}\sum_{l=0}^{n}{m+r \brace k+r}_{r,\lambda}\binom{n}{l}(k-m\lambda)_{n-l,\lambda}(y)_{k,\lambda}X^{k}\frac{1}{e_{\lambda}(x)}\big(XD+r\big)_{l,\lambda}e_{\lambda}^{y-k\lambda}(x)\nonumber\\
&=\sum_{k=0}^{m}\sum_{l=0}^{n}{m+r \brace k+r}_{r,\lambda}\binom{n}{l}(k-m\lambda)_{n-l,\lambda}(y)_{k,\lambda}x^{k}B_{l,\lambda}^{(r)}(x,y-k\lambda). \nonumber
\end{align}
Therefore, by \eqref{41} and \eqref{49}, we obtain Spivey-type recurrence relations for two variable fully degenerate $r$-Bell polynomials and fully degenerate $r$-Bell polynomials.
\begin{theorem}
For $m,n\ge 0$, we have
\begin{displaymath}
B_{n+m,\lambda}^{(r)}(x,y)=\sum_{k=0}^{m}\sum_{l=0}^{n}\binom{n}{l}{m+r \brace k+r}_{r,\lambda}(k-m\lambda)_{n-l,\lambda}(y)_{k,\lambda}x^{k}B_{l,\lambda}^{(r)}(x,y-k\lambda),
\end{displaymath}
and
\begin{displaymath}
B_{n+m,\lambda}^{(r)}(x)=\sum_{k=0}^{m}\sum_{l=0}^{n}\binom{n}{l}{m+r \brace k+r}_{r,\lambda}(k-m\lambda)_{n-l,\lambda}(1)_{k,\lambda}x^{k}B_{l,\lambda}^{(r)}(x,1-k\lambda).
\end{displaymath}
\end{theorem}
\begin{remark}
Letting $\lambda \rightarrow 0$, we obtain
\begin{equation*}
\phi_{n+m}^{(r)}(x,y)=\sum_{k=0}^{m}\sum_{l=0}^{n}\binom{n}{l}{m+r \brace k+r}_{r}k^{n-l}y^{k}x^{k}\phi_{l}^{(r)}(x,y),
\end{equation*}
and
\begin{equation*}
\phi_{n+m}^{(r)}(x)=\sum_{k=0}^{m}\sum_{l=0}^{n}\binom{n}{l}{m+r \brace k+r}_{r}k^{n-l}x^{k}\phi_{l}^{(r)}(x),
\end{equation*}
where $\phi_{n}^{(r)}(x,y)$ and $\phi_{n}^{(r)}(x)$ are the two variable $r$-Bell polynomials and the $r$-Bell polynomials, respectively given by (see \eqref{33}, \eqref{39})
\begin{displaymath}
e^{x(ye^{t}-1)}e^{rt}=\sum_{n=0}^{\infty}\phi_{n}^{(r)}(x,y)\frac{t^{n}}{n!},\quad
e^{x(e^{t}-1)}e^{rt}=\sum_{n=0}^{\infty}\phi_{n}^{(r)}(x)\frac{t^{n}}{n!}.
\end{displaymath}
\end{remark}
\section{Conclusion}
We successfully introduced new fully degenerate Bell polynomials, $B_{n,\lambda}(x)$, to address the issue of extraneous polynomials in a prior Spivey-type recurrence relation. By employing the operators $X$ and $D$, we derived natural and comprehensive Spivey-type recurrence relations for $B_{n,\lambda}(x)$ and $B_{n,\lambda}(x,y)$ (Theorem 2.8), and for their $r$-counterparts, $B_{n,\lambda}^{(r)}(x)$ and $B_{n,\lambda}^{(r)}(x,y)$  (Theorem 2.15). Furthermore, we established their fundamental properties, including Dobinski-like formulas and operator expressions. This work not only provides a corrected and more natural framework for studying the recurrence relations of degenerate Bell polynomials but also lays a foundation for future investigations into other algebraic and combinatorial properties of these new degenerate polynomial families.

\end{document}